\documentclass[reqno,a4paper]{amsart}
\usepackage{amssymb,url}
\usepackage[hypertex]{hyperref}%,pagebackref
\date{Completed on 11 June 2007}
\date{}
\allowdisplaybreaks[4]
\theoremstyle{plain}
\newtheorem{thm}{Theorem}
\newtheorem{lem}{Lemma}
\theoremstyle{remark}
\newtheorem{rem}{Remark}
\DeclareMathOperator{\td}{d\mspace{-2mu}}
\newcommand{\tn}{\mathbb{N}}

\begin{document}

\title[New proofs of the complete monotonicity of a function]
{Two new proofs of the complete monotonicity of a function involving the psi function}

\author[F. Qi]{Feng Qi}
\address[F. Qi]{Research Institute of Mathematical Inequality Theory, Henan Polytechnic University, Jiaozuo City, Henan Province, 454010, China} \email{\href{mailto: F. Qi <qifeng618@gmail.com>}{qifeng618@gmail.com}, \href{mailto: F. Qi <qifeng618@hotmail.com>}{qifeng618@hotmail.com}, \href{mailto: F. Qi <qifeng618@qq.com>}{qifeng618@qq.com}}
\urladdr{\url{http://qifeng618.spaces.live.com}}

\author[B.-N. Guo]{Bai-Ni Guo}
\address[B.-N. Guo]{School of Mathematics and Informatics,
Henan Polytechnic University, Jiaozuo City, Henan Province, 454010, China}
\email{\href{mailto: B.-N. Guo <bai.ni.guo@gmail.com>}{bai.ni.guo@gmail.com}, \href{mailto: B.-N. Guo <bai.ni.guo@hotmail.com>}{bai.ni.guo@hotmail.com}}
\urladdr{\url{http://guobaini.spaces.live.com}}

\begin{abstract}
In the present paper, we give two new proofs for the necessary and sufficient condition $\alpha\le1$ such that the function $x^\alpha[\ln x-\psi(x)]$ is completely monotonic on $(0,\infty)$.
\end{abstract}



\thanks{This paper was typeset using \AmS-\LaTeX}

\maketitle

\section{Introduction}

Recall \cite{widder} that a function $f$ is said to be completely monotonic on an interval $I$ if $f$ has derivatives of all orders on $I$ and
\begin{equation}\label{cmf-ineq}
(-1)^{n}f^{(n)}(x)\geq0
\end{equation}
for all $x\in I$ and $n\in \mathbb{N}\cup\{0\}$. The well-known Bernstein's Theorem in \cite[p.~160, Theorem~12a]{widder} states that a function $f$ on $[0,\infty)$ is completely monotonic if and only if there exists a bounded and non-decreasing function $\alpha(t)$ such that
\begin{equation}
f(x)= \int_0^\infty e^{-xt}\td\alpha(t)
\end{equation}
converges for $x\in[0,\infty)$.
\par
Recall also \cite{Atanassov, CBerg, grin-ismail, compmon2, minus-one} that a positive function $f$ is said to be logarithmically completely monotonic on an interval $I$ if $f$ has derivatives  of all orders on $I$ and
\begin{equation}\label{lcmf-ineq}
(-1)^n[\ln f(x)]^{(n)}\ge 0
\end{equation}
for all $ x\in I$ and $n\in\mathbb{N}$.
\par
It was proved explicitly in \cite{CBerg, compmon2, minus-one} and other articles that a logarithmically completely monotonic function must be completely monotonic. For more information on the logarithmically completely monotonic functions, please refer to \cite{CBerg, e-gam-rat-comp-mon} and related references therein.
\par
It is well-known that the Euler gamma function is defined by
\begin{equation}\label{egamma}
\Gamma(z)=\int^\infty_0t^{z-1} e^{-t}\td t
\end{equation}
for $\Re z>0$. The logarithmic derivative of $\Gamma(z)$, denoted by $\psi(z)=\frac{\Gamma'(z)}{\Gamma(z)}$, is called the psi or digamma function, and $\psi^{(k)}$ for $k\in\mathbb{N}$ are called the polygamma functions.
\par
In \cite{Anderson}, the function
\begin{equation}\label{theta-dfn-psi}
\theta(x)=x[\ln x-\psi(x)]
\end{equation}
was proved to be decreasing and convex on $(0,\infty)$, with two limits
\begin{equation}\label{2limitss}
\lim_{x\to0^+}\theta(x)=1\quad \text{and}\quad \lim_{x\to\infty}\theta(x)=\frac12
\end{equation}
were presented complicatedly.
\par
In \cite[p.~374]{psi-alzer}, it was pointed out that the limits in~\eqref{2limitss} can follow immediately from the representations
\begin{equation*}
\theta(x)=x\ln x-x\psi(x+1)+1\quad\text{and}\quad
\theta(x)=\frac12+\frac1{12x}-\frac{\tau}{120x^3}
\end{equation*}
for $x>0$ and $\tau\in(0,1)$.
\par
From~\eqref{2limitss} and the decreasing monotonicity of $\theta(x)$, the inequality
\begin{equation}\label{lnx-psi}
\frac1{2x}<\ln x-\psi(x)<\frac1x
\end{equation}
for $x>0$ is concluded. This extends a result in \cite{minc}, which says that the inequality~\eqref{lnx-psi} is valid for $x>1$. Refinements and generalizations of~\eqref{lnx-psi} were given in \cite{gordon, sandor-gamma-2, poly-comp-mon.tex} and related references therein. For more information, please refer to \cite{bounds-two-gammas.tex} and related references therein.
\par
In \cite{note-on-li-chen.tex}, by employing the monotonicity of $\theta(x)$, it was recovered simply that the double inequality
\begin{equation}\label{li-chen-ineq}
\frac{x^{x-\gamma}}{e^{x-1}}<\Gamma(x)<\frac{x^{x-1/2}}{e^{x-1}}
\end{equation}
holds for $x>1$, the constants $\gamma$ and $\frac12$ are the best possible, the left-hand side inequality in~\eqref{li-chen-ineq} holds also for $0<x<1$, but the right-hand side inequality in~\eqref{li-chen-ineq} reverses, where $\gamma$ is Euler-Mascheroni's constant. Furthermore, by virtue of the decreasing monotonicity and convexity of $\theta(x)$, it was showed in \cite{note-on-li-chen.tex} that the function
\begin{equation}
h(x)=\frac{e^x\Gamma(x)}{x^{x-\theta(x)}}
\end{equation}
on $(0,\infty)$ has a unique maximum $e$ at $x=1$, with two limits
\begin{equation}\label{2limits}
\begin{aligned}
\lim_{x\to0^+}h(x)&=1 & \text{and} &&\lim_{x\to\infty}h(x)&=\sqrt{2\pi}\,.
\end{aligned}
\end{equation}
Consequently, three sharp inequalities
\begin{equation}\label{qi-guo-zhang-ineq-1}
\frac{x^{x-\theta(x)}}{e^x}<{\Gamma(x)}\le \frac{x^{x-\theta(x)}}{e^{x-1}}
\end{equation}
on $(0,1]$,
\begin{equation}\label{qi-guo-zhang-ineq}
\frac{\sqrt{2\pi}\,{x^{x-\theta(x)}}}{e^x}
<{\Gamma(x)}\le\frac{x^{x-\theta(x)}}{e^{x-1}}
\end{equation}
on $[1,\infty)$, and
\begin{equation}\label{note-li-chen-ineq}
I(x,y)<\biggl\{\frac{x^{\theta(x)}\Gamma(x)} {y^{\theta(y)}\Gamma(y)}\biggr\}^{1/(x-y)}
\end{equation}
for $x\ge1$ and $y\ge1$ with $x\ne y$, where
\begin{equation}
I(a,b)=\frac1e\biggl(\frac{b^b}{a^a}\biggr)^{1/(b-a)}
\end{equation}
for $a>0$ and $b>0$ with $a\ne b$ is called the identric or exponential mean, are deduced directly. If $0<x\le1$ and $0<y\le1$ with $x\ne y$, the inequality~\eqref{note-li-chen-ineq} is reversed.
\par
In \cite[pp.~374--375, Theorem~1]{psi-alzer}, by using the well-known Binet's formula and complicated calculating techniques for integrals, the monotonicity and convexity of $\theta(x)$ was extended to the complete monotonicity: For real number $\alpha$, the function
\begin{equation}\label{theta-alpha}
\theta_\alpha(x)=x^\alpha[\ln x-\psi(x)]
\end{equation}
is completely monotonic on $(0,\infty)$ if and only if $\alpha\le1$.
\par
The aim of this paper is to give two new proofs of the complete monotonicity of the function $\theta_\alpha(x)$, which can be restated as the following Theorem~\ref{theta-comp-mon}, since this function $\theta_\alpha(x)$ has many meaningful applications as stated above.

\begin{thm}\label{theta-comp-mon}
For real number $\alpha$, the function $\theta_\alpha(x)$ defined by~\eqref{theta-alpha} is completely monotonic on $(0,\infty)$ if and only if $\alpha\le1$, with two limits
\begin{equation}\label{a=1-limits}
\lim_{x\to0^+}\theta_1(x)=1,\quad \lim_{x\to\infty}\theta_1(x)=\frac12
\end{equation}
and, for $\alpha<1$,
\begin{equation}\label{limits-2}
\lim_{x\to0^+}\theta_\alpha(x)=\infty,\quad \lim_{x\to\infty}\theta_\alpha(x)=0.
\end{equation}
\end{thm}

\section{Lemmas}

In order to prove Theorem~\ref{theta-comp-mon}, the following lemmas are needed.

\begin{lem}[\cite{abram}]
For $i\in\mathbb{N}$, $x>0$, $a>0$ and $b>0$,
\begin{align}\label{psisymp4}
\psi^{(i-1)}(x+1)&=\psi^{(i-1)}(x)+\frac{(-1)^{i-1}(i-1)!}{x^i},\\
\ln\frac{b}a&=\int_0^\infty\frac{e^{-at}-e^{-bt}}t\td t,
\label{abram-230-5.1.32} \\
\psi ^{(i)}(x)&=(-1)^{i+1}\int_{0}^{\infty}\frac{t^{i} e^{-xt}}{1-e^{-t}}\td t, \label{psim}\\
\label{binetcor}
\psi(x)-\ln x+\frac1x&=\int_{0}^{\infty}\left(\frac{1}{t}-\frac{1}{e^t-1}\right)e^{-xt}\td t.
\end{align}
\end{lem}

\begin{lem}[\cite{poly-comp-mon.tex}]
For $x>0$,
\begin{gather}\label{psi(x+1)ineq}
\frac{1}{2x}-\frac{1}{12x^2}<\psi(x+1)-\ln x<\frac{1}{2x},\\
\label{lem2}
\frac{1}{2x^2}-\frac{1}{6x^3}<\frac{1}{x}-\psi '(x+1)
<\frac{1}{2x^2}-\frac{1}{6x^3}+\frac{1}{30x^5}.
\end{gather}
\end{lem}

\begin{lem}\label{comp-thm-1}
Inequalities
\begin{equation}
\ln x-\frac1x\le\psi(x)\le\ln x-\frac1{2x}
\end{equation}
and
\begin{equation}\label{qi-psi-ineq}
\frac{(k-1)!}{x^k}+\frac{k!}{2x^{k+1}}\le
(-1)^{k+1}\psi^{(k)}(x)\le\frac{(k-1)!}{x^k}+\frac{k!}{x^{k+1}}
\end{equation}
hold on $(0,\infty)$ for $k\in\tn$.
\end{lem}

\begin{proof}
In \cite{sandor-gamma-2}, the function $\psi(x)-\ln x+\frac{\alpha}x$ was proved to be completely monotonic on $(0,\infty)$ if and only if $\alpha\ge1$ and so was its negative if and only if $\alpha\le\frac12$. In \cite{chen-qi-log-jmaa}, the function $\frac{e^x\Gamma(x)} {x^{x-\alpha}}$ was proved to be logarithmically completely monotonic on $(0,\infty)$ if and only if $\alpha\ge1$ and so was its reciprocal if and only if $\alpha\le\frac12$. From these, considering~\eqref{cmf-ineq} and~\eqref{lcmf-ineq}, inequalities in~\eqref{qi-psi-ineq} are derived straightforwardly.
\end{proof}

\begin{lem}\label{elbert-laforgia-lem}
If $f(x)$ is a function defined in an infinite interval $I$ such that
$f(x)-f(x+\varepsilon)>0$ and $\lim_{x\to\infty}f(x)=\delta$ for $x\in I$ and
some $\varepsilon>0$, then $f(x)>\delta$ in $I$.
\end{lem}

\begin{proof}
By induction, for any $x\in I$,
\begin{equation*}
f(x)>f(x+\varepsilon)>f(x+2\varepsilon)>\dotsm>f(x+k\varepsilon)\to\delta
\end{equation*}
as $k\to\infty$. The proof of Lemma \ref{elbert-laforgia-lem} is complete.
\end{proof}

\begin{rem}
Lemma~\ref{elbert-laforgia-lem} is simple, but it is very effectual in dealing with some problems concerning (logarithmically) completely monotonic properties of functions involving the gamma, psi, polygamma functions.
\end{rem}

\section{The first proof of Theorem~\ref{theta-comp-mon}}

Straightforward computation gives
\begin{align*}
\theta_1(x+1)-\theta_1(x)&=(x+1)\ln(x+1)-x\ln x+x[\psi(x)-\psi(x+1)]-\psi(x+1)\\
&=(x+1)\ln(x+1)-x\ln x-\psi(x+1)-1
\end{align*}
and
\begin{align*}
[\theta_1(x+1)-\theta_1(x)]'&=\ln(x+1)-\ln x-\psi'(x+1)\\
&=\int_0^\infty\biggl[\frac{1-e^{-t}}t -\frac{te^{-t}}{1-e^{-t}}\biggr]e^{-xt}\td t\\
&=\int_0^\infty\frac{e^{-t}+e^t-t^2-2}{t(e^t-1)}e^{-xt}\td t\\
&>0
\end{align*}
by using formulas~\eqref{psisymp4}, \eqref{abram-230-5.1.32} and~\eqref{psim}. Hence,
\begin{equation*}
(-1)^i[\theta_1(x+1)-\theta_1(x)]^{(i)}
=\bigl[(-1)^i\theta_1^{(i)}(x+1)]-[(-1)^i\theta_1^{(i)}(x)\bigr]<0
\end{equation*}
on $(0,\infty)$ for $i\in\tn$.
\par
Using the inequality~\eqref{psi(x+1)ineq} yields
\begin{multline*}
(x+1)\ln\biggl(1+\frac1x\biggr)-\frac1{2x}-1<\theta_1(x+1)-\theta_1(x)\\*
<(x+1)\ln\biggl(1+\frac1x\biggr)-\frac1{2x}+\frac1{12x^2}-1,
\end{multline*}
which implies that $\lim_{x\to\infty}[\theta_1(x+1)-\theta_1(x)]=0$. Since the
function $\theta_1(x+1)-\theta_1(x)$ is increasing on $(0,\infty)$, it is obtained
that $\theta_1(x+1)-\theta_1(x)<0$ on $(0,\infty)$.
\par
Utilizing~\eqref{psisymp4} and~\eqref{psi(x+1)ineq} leads easily to
$\lim_{x\to\infty}\theta_1(x)=\frac12$.
\par
Direct calculation gives $\theta_1'(x)=\ln x-\psi(x)-x\psi'(x)+1$ and
$$
\theta_1^{(i)}(x)=\frac{(-1)^i(i-2)!}{x^{i-1}}-i\psi^{(i-1)}(x)-x\psi^{(i)}(x)
$$
for $i\ge2$. Combination of~\eqref{psisymp4}, \eqref{psi(x+1)ineq} and~\eqref{lem2} yields that $\lim_{x\to\infty}\theta_1'(x)=0$. The inequality~\eqref{qi-psi-ineq} means that $\lim_{x\to\infty}\theta_1^{(i)}(x)=0$ for $i\ge2$.
\par
By the above argument and Lemma~\ref{elbert-laforgia-lem}, it is concluded that
$(-1)^k\theta_1^{(k)}(x)\ge0$ on $(0,\infty)$ for $k\ge0$, which means that the
function $\theta_1(x)$ is completely monotonic on $(0,\infty)$ with
$\lim_{x\to\infty}\theta_1(x)=\frac12$.
\par
The validity of the limit $\lim_{x\to0^+}\theta_1(x)=1$ follows from the formula~\eqref{binetcor}.
\par
It is clear that $\theta_\alpha(x)=x^{\alpha-1}\theta_1(x)$ and $x^{\alpha-1}$ is also completely monotonic on $(0,\infty)$ for $\alpha<1$. Since the product of any finite completely monotonic functions on an interval $I$ is also completely monotonic on $I$, the function $\theta_\alpha(x)$ is completely monotonic on $(0,\infty)$ for $\alpha<1$.
\par
Conversely, if the function $\theta_\alpha(x)$ is completely monotonic on $(0,\infty)$, then $\theta_\alpha(x)$ is decreasing and positive on $(0,\infty)$. From the formula~\eqref{psisymp4} and the inequality~\eqref{psi(x+1)ineq}, it follows that
\begin{equation}\label{psi(x+1)ineq-var}
\frac{1}{2x}+\frac{1}{12x^2}>\ln x-\psi(x)>\frac{1}{2x}
\end{equation}
and
\begin{equation}\label{lnx-psix}
\frac{1}{2x^{1-\alpha}}+\frac{1}{12x^{2-\alpha}}>x^\alpha[\ln x-\psi(x)]>\frac{1}{2x^{1-\alpha}}
\end{equation}
for $x>0$, which means that $x^\alpha[\ln x-\psi(x)]$ tends to $\infty$ as $x\to\infty$ if $\alpha>1$. This contradicts with the decreasingly monotonic property of $\theta_\alpha(x)$ on $(0,\infty)$. Hence, the necessary condition $\alpha\le1$ follows.
\par
It is obvious that the inequality~\eqref{lnx-psix} implies the two limits in~\eqref{limits-2}. The proof of Theorem~\ref{theta-comp-mon} is complete.

\section{The second proof of Theorem~\ref{theta-comp-mon}}

Let
\begin{equation}
h(t)=\frac{1}{t}-\frac{1}{e^t-1}=\frac{e^t-1-t}{t(e^t-1)}
\end{equation}
for $t\ne0$ and $h(0)=\frac12$. Integration by part in~\eqref{binetcor} yields
\begin{equation}
\begin{split}\label{psi(x)-ln x+frac1x}
\psi(x)-\ln x+\frac1x&=-\frac1x\biggl\{\bigl[h(t)e^{-xt}\bigr]\big|_{t=0}^{t=\infty} -\int_{0}^{\infty}h'(t)e^{-xt}\td t\biggr\}\\
&=\frac1{2x}+\frac1x\int_{0}^{\infty}h'(t)e^{-xt}\td t.
\end{split}
\end{equation}
Multiplying on all sides of~\eqref{psi(x)-ln x+frac1x} by $x$ and rearranging gives
\begin{equation}
x[\ln x-\psi(x)]=-\frac12+\int_{0}^{\infty}h'(t)e^{-xt}\td t.
\end{equation}
In~\cite{best-constant-one-simple.tex, best-constant-one.tex, best-constant-one-simple-real-AMEN.tex} and related references therein, the function $h(t)$ was shown to be decreasing on $(-\infty,\infty)$, concave on $(-\infty,0)$ and convex on $(0,\infty)$. This means that the function $\theta_1(x)$ is completely monotonic on $(0,\infty)$ and that the second limit in~\eqref{a=1-limits} follows.
This means that if $\alpha>1$ then the function $\theta_\alpha(x)=x^{\alpha-1}\theta_1(x)$ tends to infinity for $x$ tending to infinity and therefore it cannot be completely monotonic, that is, the condition $\alpha\le1$ is necessary. The second proof of Theorem~\ref{theta-comp-mon} is complete.

\begin{rem}
The second proof of Theorem~\ref{theta-comp-mon} can also be demonstrated as follow.
It is easy to see that
\begin{equation}
\frac1x=\int_0^\infty e^{-xu}\td u,\quad x>0.
\end{equation}
Substituting it into~\eqref{binetcor} gives
\begin{equation}
\ln x-\psi(x)=\int_0^\infty\biggl(\frac1{1-e^{-t}}-\frac1t\biggr)e^{-xt}\td t \triangleq\int_0^\infty\rho(t)e^{-xt}\td t.
\end{equation}
An integration by part and a multiplication by $x$ yields
\begin{equation}
x[\ln x-\psi(x)]=\frac12+\int_0^\infty\rho'(t)e^{-xt}\td t,
\end{equation}
where
\begin{equation}
\rho'(t)=\frac{1}{t^2}-\frac{e^{-t}}{(1-e^{-t})^2}
=\frac{2e^{-t}}{t^2(1-e^{-t})^2}\biggl(\frac{e^t+e^{-t}}2-1-\frac{t^2}2\biggr).
\end{equation}
Making use of the power series expansion of $e^t$ at $t=0$ reveals easily that $\rho'(t)$ is positive on $(0,\infty)$. So the function $\theta_1(x)$ is completely monotonic with the limit $\frac12$ at infinity.
\end{rem}

\end{document}